\input amstex
\documentstyle{amsppt}
\nologo
\document
\loadmsam
\magnification 1200
\pageheight{7.25 in}

\define \modq{\mod q}

\define\zp{\Bbb Z_p}

\define\bbbn{\Bbb N}

\define\r {\Cal R}
\define\la {\lambda}
\define\p{\prime}
\define\kl{\circ}
\topmatter
\title Some Hereditarily Just Infinite Subgroups of $\Cal J$ \endtitle
\author Cornelius Griffin \endauthor
\affil School of Mathematical Sciences \\ University of Nottingham \\ 
Nottingham \\ NG7 2RQ \endaffil
\date{10th September 2002} \enddate
\abstract
This work examines the commutator structure of some closed subgroups of the 
wild group of automorphisms of a local field with perfect residue field, a 
group we call $\Cal J.$ In particular, we establish a new approach to 
evaluating commutators in $\Cal J$ and using this method investigate the 
normal subgroup structure of some classes of index subgroups of $\Cal J$ as 
introduced by Klopsch. We provide new proofs of Fesenko\rq s results that 
lead to a proof that the torsion free group $T =\{ t+\sum_{k\geq 1} 
a_kt^{qk+1}: a_k \in \Bbb F_p\}$ is hereditarily just infinite, and by 
extending his work, we also demonstrate the existence of a new class of 
hereditarily just infinite subgroups of $\Cal J$ which have non-trivial 
torsion.
\endabstract
\endtopmatter

\leftheadtext{Cornelius Griffin}
\rightheadtext{Index Subgroups of $\Cal J$}

\heading {Index Subgroups of  $\Cal J$}
\endheading

\heading{\bf {Section 1: Introduction}}
\endheading
\vskip 0.25in

Although as a result of the work of Camina \cite{C1} it is known that $\Cal
J$ has an extremely rich subgroup structure, it was generally considered
difficult to write down arbitrary subgroups of $\Cal J.$ This problem
recieved some attention from Klopsch, who in his thesis \cite{K} considered
a particular class of subgroups in $\Cal J$ which are easier to describe
than others. His class was defined as follows: take a subset $\Lambda
\subset \bbbn$ and look at the equivalent subset

$$
\Cal J (\Lambda) := \left\{t+\sum_{\la \in \Lambda} a_{\la}t^{\la +1}:
a_{\la}
\in \Bbb F_p\right\}
$$
of $\Cal J.$ Under certain arithmetic conditions on $\Lambda$ it was already
known that this subset actually formed a subgroup of $\Cal J.$ Klopsch
managed to give a complete description of sets $\Lambda$ for which $\Cal
J(\Lambda)$ is a subgroup of $\Cal J.$ He proved

\proclaim {Theorem \cite{K}}
Let $\Lambda$ and $\Cal J(\Lambda)$ be as above. Then $\Cal J(\Lambda) \leq
\Cal J$ if and only if

$$
\forall \la, \mu \in \Lambda, \forall k \leq \la +1: \binom {\la+1}k \equiv
   0   \mod p \text { or } \la +k\mu \in \Lambda.
$$
\endproclaim

Given this arithmetic condition he went on to describe four classes of
subsets which are easily demonstrable to satisfy the condition in the
Theorem:

\roster
\item A: = $d\bbbn$ for some $d \in \bbbn;$
\item B: = $p\bbbn \cup p\bbbn -1;$
\item C: = $p^i\bbbn -1$ for some $i \in \bbbn;$
\item D: = $ \{p^i-1: i \in \bbbn \}.$
\endroster

Klopsch went on to prove some elementary results about the groups arising
from these subsets of $\bbbn$. He considered properties such as finite
generation,
centralizers, normalizers etc., and also used these subgroups to calculate
part of the Hausdorff Spectrum for the Nottingham Group.

Independently of Klopsch\rq s work, Fesenko has studied the group he called
$T:=\Cal J (q\bbbn)$ for some power $q=p^r$ of $p.$ It transpires that this
group has some interesting properties, both as a group in its own right, and
also as the Galois group of a particular type of field extension. Coates and
Greenberg had asked the following

\proclaim{Question \cite{CG}}
Is it true that for every finite extension $K$ of $\Bbb Q_p$ there exists a
deeply ramified Galois $p$-extension $M$ of $K$ so that no subfield $M^{\p}$
of $M$ is an infinitely ramified Galois extension of a finite extension $Q$
of $\Bbb Q_p$ with $Gal(M^{\p}/Q)$ being a $p$-adic Lie group?
\endproclaim

The terms need not concern us greatly here; I would just point out that an
arithmetically profinite field extension is a particular type of a deeply 
ramified field extension. In particular an arithmetically profinite field 
extension is one
whose non-trivial upper ramification groups are all open --- this is the
property we need. Fesenko managed to give an affirmative answer to this
question by taking a field extension whose Galois Group is equal to $T.$
Then to answer the question, the problem is reduced to a need to demonstrate
some properties of the group $T.$ Fesenko proved

\proclaim {Theorem \cite{F1}}
Let $q$ be a fixed power of the odd prime $p$ and let $T$ be the group
$$
\left\{ t+\sum_{k=1}^{\infty} a_kt^{qk+1} : a_k \in \Bbb F_p \quad \forall k
\right\}.
$$
Then
\roster
\item If $\sigma \in T_i\backslash T_{i+1}$ then $\sigma^p \in
T_{pi}\backslash T_{pi+1}$ and so $T$ is torsion-free$;$
\item $[T_i,T_i] \leq T_{(q+1)i+1}$ and the group $T_i /T_{i+1}$ is abelian
of exponent $p;$
\item $[T,T]T^p >T_{q+2};$ the number of generators of $T$ is at most $q+1;$
\item $T$ is not $p$-adic analytic$;$
\item $T$ is hereditarily just infinite$;$ i.e., non-trivial closed normal
subgroups of open subgroups are open.
\endroster
\endproclaim

By using this Theorem and some ramification theory it is possible to give an
answer to the Coates and Greenberg Question. For instance, that $T$ is
hereditarily just infinite allows one to show that the corresponding field
extension must be arithmetically profinite.

This Theorem of Fesenko can be seen in other contexts as well. Hereditarily
just infinite groups, and their ancestors the just infinite groups play the
same role in pro-$p$ group theory as the simple groups play in finite group
theory. So it is natural to try to classify the hereditarily just infinite
groups as one does for the finite simple groups. In the advent of Fesenko\rq
s work, the current state of this classification is as follows:

\proclaim {Partition of JI Groups}
The class of just infinite pro-$p$ groups consists of four types:
\roster
\item Solvable ones that are linear over $\zp;$
\item Nonsolvable ones that are linear over $\zp;$
\item Nonsolvable ones that are linear over $\Bbb F_p[[t]];$
\item The rest --- namely groups of  \lq\lq Nottingham type\rq\rq , \lq\lq
Grigorchuk type\rq\rq , \lq\lq Fesenko type\rq\rq , and so on.
\endroster
\endproclaim

Apart from aesthetic concerns there are strong practical reasons why it
would be useful to complete this classification. For instance, Boston has
reduced the Fontaine-Mazur Conjecture concerning $p$-adic Galois
representations to a question about hereditarily just infinite pro-$p$
groups \cite{B}. As the Fontaine-Mazur Conjecture implies amongst other
things an alternative proof of Fermat\rq s Theorem, one can see the strength
of such a classification of HJI groups.

To prove that a group is HJI, the most important thing a person must do is
to understand properly the commutator structure of the group in question.
The main thrust of Fesenko\rq s argument is a series of highly technical
Lemmas in which sufficiently many commutators are evaluated to leave him
able to deduce his final results.

The aim of this paper is twofold. Firstly, we introduce a new approach to
calculating commutators in groups of formal power series. To illustrate the
usefulness of this approach, we provide new proofs of the Lemmas contained
in the work of Fesenko about commutators in $T$ and in the interest of
completeness, we indicate how from here one may deduce that $T$ is HJI.

Secondly, we use the same methods, and indeed Fesenko\rq s proof that $T$ is 
HJI to evaluate commutators in $\Cal J (B)$
and as a result, we are able to prove the following

\proclaim{Theorem 6.3} Let $B:= p\bbbn \cup p\bbbn -1$ and let $S:=\Cal
J(B).$ Then $S$ is a hereditarily just infinite pro-$p$ group with 
non-trivial torsion.
\endproclaim

The method employed could also be used to study the normal subgroup
structure of subgroups of the Nottingham Group in more generality. It may be
possible to extend the methods used here to calculate, for instance, the
obliquity/width of various index subgroups. This has been done for the group
$S$ in a recent preprint of Barnea and Klopsch \cite{BK} which also contains
an alternative proof of Theorem 6.3.

This work is part of the author\rq s PhD thesis, Nottingham 2002, supported 
by an EPSRC studentship and under the supervision of Professor Ivan Fesenko. 
The author thanks EPSRC and also B. Klopsch, I. Fesenko for their interest 
and helpful comments.

\vskip 0.5in
\heading {\bf {Section 2: Evaluating Commutators in Groups of Formal Power
Series}}
\endheading
\vskip 0.25in

Throughout this text, by $[v,u]$ we will mean $v \kl u\kl v^{-1}\kl u^{-1}.$
Here $\kl$ denotes the group operation; to indicate formal products of two
power series $v$ and $u$ we will merely write $vu.$
\newline

So suppose we are given two formal power series $u$ and $v$ in some group of
formal power series. Then we may write
$$
[v,u] := t+\sum_{k \geq l}a_kt^{k+1}
$$
for some $a_k \in \Bbb F_p$ and some $l \in \bbbn.$ We want to understand
the series $(a_k)$ for values $k \in \bbbn.$ Now $[v,u]=v\kl u\kl v^{-1}\kl
u^{-1}$ and so
$$
v\kl u\kl v^{-1}\kl u^{-1} = t+\sum_{k\geq l}a_kt^{k+1}
$$
$$
\iff \quad v\kl u\kl v^{-1}\kl u^{-1}\kl u\kl v = (t+\sum_{k\geq
l}a_kt^{k+1}) \kl u\kl v
$$
$$
\iff \quad \quad v\kl u -u\kl v = \sum_{k\geq l} a_k(u \kl v)^{k+1}.
$$
So we have reduced understanding commutators of elements to the simpler
problem of understanding products of elements. Now to evaluate commutators
of formal power series, we need only solve some recurrence relations on the
$a_k.$ To illustrate this process I will now give a new proof of the Lemmas
used by Fesenko evaluating  particular commutators that were necessary to
show  that $T$ is HJI. In all the work that follows, we will implicitly use
the following two Lemmas.

\proclaim{Crucial Lemma 2.1} Let $m \geq n \in \bbbn$, let $l_1, \dots
,l_n \in \bbbn$ be so that $\sum_{i=1}^n l_i =m.$ Write $A_{m,n}(l_1, \dots
,l_n)$ for the number of maps

$$
f:\{1, \dots ,m\} \to \{1, \dots ,n\} \text { so that } |f^{-1}(i)|=l_i.
$$
Then

$$
A_{m,n}(l_1, \dots ,l_n) = \prod_{i=1}^n {\binom
{m-\sum_{k=1}^{i-1}l_k}{l_i}}.
$$
\endproclaim

\demo{Proof}
When $n=1$ the result is obvious.
\newline
Now count the number of maps $f:\{1, \dots ,m\} \to \{1, \dots ,n+1\} \text
{ so that } |f^{-1}(i)|=l_i.$ To do this, take $l_1$ elements in $\{1, \dots
,m\}$ that map onto $1.$ This can be done in $\binom {m}{l_1}$ ways and so
the number of maps having the required property is $\binom {m}{l_1}$ times
the number of maps $f:\{1, \dots, m-l_1\} \to \{1 , \dots , n\}$ so that
$|f^{-1}(i)|=l_i$ where we have done some relabelling. The result now
follows by induction on $n.$
\enddemo

\proclaim{Lemma 2.2 \cite{Lu}} Let $n \in \bbbn$ and let $a=\sum_{i=0}^n
a_ip^i, b=\sum_{i=0}^nb_ip^i \in \bbbn,$ where $a_i,b_i \in \{0,\dots
,p-1\}$ for all $i\in \{0,\dots ,p-1\}.$ Then
$$
\binom ab \not\equiv   0   \mod p \iff a_i\geq b_i \quad \forall \, i.
$$
\endproclaim

Throughout the following work we will use the following notation; for a
given element $u =t+\sum_{k\geq 1} u_{k+1}t^{k+1}$ we will set
$$
E_i(u) := u_i \text { and }
S(u) := \{ i \in \bbbn : E_i(u) \not= 0 \}.
$$

\vskip 0.5in
\heading {\bf {Section 3: Commutators in $T$}}
\endheading
\vskip 0.25in

$T$ is the group consisting of formal power series of the form
$$
t+\sum_{k\geq 1}a_{qk+1}t^{qk+1}, \quad a_l \in \Bbb F_p
$$
for $q$ some power of the prime $p.$ We can define a filtration on $T$ by
setting
$$
T_i:=\left\{t+\sum_{k\geq i}a_{qk+1}t^{qk+1} \right\}.
$$
In order to investigate the commutator structure of this group we want to
calculate $[v,u]$ for arbitrary elements $v \in T_j,u \in T_i$ for $i \geq
j.$
\newline
For convenience set $i=j+e$ and given elements $u,v \in T$ write $(u\kl
v)_{nc}=u\kl v -u-v+t.$ Notice this means that $(v\kl u)_{nc} -(u\kl
v)_{nc}=v\kl u -u\kl v.$
\newline

Recall the reasoning outlined that will allow us to calculate commutators of
this form. We have that  $[v,u] = t+ \sum_{k \geq K} a_{qk+1}t^{qk+1}$ if
and only if
$$
v \kl u-u\kl v = \sum_{k\geq K}a_{qk+1}(u\kl v)^{qk+1}
$$
and so we must evaluate the compositions $v\kl u, u\kl v.$ It is a simple
exercise to verify that we may write
$$
u\kl v=u+v-t+\sum_{s\geq2} {\left(\sum_{k=j}^{s-j}
u_{q(k+e)+1}f_{s,k}\right )t^{q(s+e)+1}}
$$
and
$$
v\kl u=v+u-t+\sum_{s\geq2} {\left (\sum_{k=j}^{s-j} v_{qk+1}g_{s,k}\right
)t^{q(s+e)+1}}
$$
where now
$$
f_{s,k}=  \sum_{j(1)+\dots +j(q(k+e)+1)=q(s+e)+1} v_{j(1)}\dots v_{j(qk+1)}
$$
and
$$
g_{s,k}=\sum_{j(1)+\dots +j(qk+1)=q(s+e)+1} u_{j(1)}\dots u_{j(qk+1)}.
$$
In  the first sum we sum over integers $j(l)$ so that $j(l)>1 \Rightarrow
j(l)\geq qj+1$  whereas the  second sum is over integers $j(l)$ so that
$j(l)>1 \Rightarrow j(l)\geq q(j+e)+1.$ It should be clear now why we have
identified  Lemmas 2.1, 2.2 as crucial if we want to evaluate these
products, and hence also commutators in $T$. The parts of the sums above
requiring analysis are the $f_{s,k},g_{s,k}.$ Let us first consider
$f_{s,k};$ we collect information into a

\proclaim{Lemma 3.1}
Fix $s$ and $k$ and take notation as above.
\roster
\item Suppose there exists a unique $l$ so that $j(l) >1.$ Set $j(l)=qm+1$
and $j(n)=1$ for all $n \not= l.$ Then $m+k=s$ and so $s \geq 2j;$
\item Suppose there exists $l_1, \dots ,l_q$ so that $j(l_i)>1.$ Set
$j(l_i)=qm+1.$ Then $qm+k=s$ and so $s \geq qj+j;$
\item Suppose there exists $l_1, \dots ,l_{q+1}$ so that $j(l_i)>1.$ Set
$j(l_i)=qm+1$ for $i=1, \dots, q$ and $j(l_{q+1}) =qn+1.$ Then $qm+n+k =s$
and so $s \geq qj+2j;$
\item Suppose in general there exist $b_1$ so that $j(l) =qm_1+1, b_2 $ so
that $j(l)=qm_2+1$ and so on to $b_d$ so that $j(l)=qm_d+1.$ Then $s=b_1m_1
+\dots +b_dm_d +k$ and so $s \geq \Delta (b_1+\dots +b_d)j +j$ where $\Delta
:=|\{i:m_i>0\}|.$
\endroster

\endproclaim

This simple Lemma allows one to write down $f_{s,k}$ for small values of $s$
and so to write down $u \kl v.$

\proclaim{Proposition 3.2}
\roster
\item Let $2j \leq s<qj+j.$ Then $f_{s,k} = v_{q(s-k)+1}$ and so the
coefficient of $t^{q(s+e)+1}$ in $(u\kl v)_{nc}$ is
$$
\sum_{k=j}^{s-j}u_{q(k+e)+1}v_{q(s-k)+1};
$$
\item Let $qj+j \leq s<qj+2j.$ Then the coefficient of  $t^{q(s+e)+1}$ in
$(u\kl v)_{nc}$ is
$$
\sum_{k=j}^{s-j} u_{q(k+e)+1}v_{q(s-k)+1} + \sum \Sb m,k\geq j, \\ qm+k=s
\endSb {(k+e)u_{q(k+e)+1}v_{qm+1}};
$$
\item Let $s=qj+2j.$ In addition to the terms described above we also get in
$(u\kl v)_{nc}$ the term
$$
(j+e)u_{q(j+e)+1}v_{qj+1}^2;
$$
\item This process continues to expand in a uniform way as $s$ increases in
size.
\endroster

\endproclaim

\demo{Proof}
Everything follows straightforwardly from Lemmas 2.1, 2.2, 3,1. To
illustrate the process I will evaluate $f_{qj+j,j}.$ We have
$$
f_{qj+j,j} = \sum_{j(1)+\dots +j(q(j+e)+1)=q(qj+j+e)+1} v_{j(1)}\dots
v_{j(qj+1)}.
$$
If there exists a unique $l$ so that $j(l) >1$ then this $j(l)$ can be
chosen in $\binom {q(j+e)+1}1$ ways and we have $j(l) +q(j+e)=q(qj+j+e)+1$
from which it follows that $j(l) =q(qj)+1$ and we get a contribution to
$f_{qj+j,j}$ of $\binom {q(j+e)+1}1 v_{q(qj)+1}.$

If there exists $d$ say values of $l$ so that $j(l)>1$ then these $d$ values
can be chosen in $\binom {q(j+e)+1}d$ ways. Thus Lucas\rq\   Lemma 2.2
tells us that $d$ is either $0$ or $1   \modq .$ So the size of $s$ tells
us that we must have $d=q$ from which we may deduce that $qj(l) +q(j+e)+1-q
= q(qj+j+e)+1.$ Thus $j(l)=qj+1$ and we get a contribution to $f_{qj+j,j}$
of $\binom {q(j+e)+1}q v_{qj+1}.$

The result follows.

\enddemo

We can do exactly the same for $v \kl u;$ omitting  the details,
one may prove

\proclaim{Proposition 3.3}
\roster
\item Let $2j \leq s<qi+j-e.$ Then $g_{s,k} = u_{q(s+e-k)+1}$ and so the
coefficient of $t^{q(s+e)+1}$ in $(v\kl u)_{nc}$ is
$$
\sum_{k=j}^{s-j}v_{q(k)+1}u_{q(s+e-k)+1};
$$
\item Let $qi+j-e \leq s<qi+2j.$ Then the coefficient of  $t^{q(s+e)+1}$ in
$(v\kl u)_{nc}$ is
$$
\sum_{k=j}^{s-j} v_{q(k)+1}u_{q(s+e-k)+1} + \sum \Sb k\geq j, \\l\geq j+e \\
ql+k=s+e \endSb {kv_{q(k)+1}u_{ql+1}};
$$
\item Let $s=qi+2j.$ In addition to the terms described above we also get in
$(v\kl u)_{nc}$ the term
$$
ju_{q(j+e)+1}^2v_{qj+1};
$$
\item This process continues to expand in a uniform way as $s$ increases in
size.
\endroster
\endproclaim

So we now have all the tools we will require to evaluate commutators in $T.$
The first thing to notice is that the leading coefficients of $[v,u]$ and
$t+v\kl u-u\kl v$ are the same and so we may immediately deduce

\proclaim{Proposition 3.4} Let $u,v$ be as above. Then
\roster
\item  $[v,u]=t-iu_{qi+1}v_{qj+1}t^{q(qj+i)+1}+\dots ;$
\item $[T_j,T_i] \leq T_{qj+i}$ and if $p$ divides $i$ then furthermore
$[T_j,T_i] \leq T_{qj+i+1}.$
\endroster
\endproclaim

\demo{Proof}
We have
$$
\split
t+v\kl u-u\kl v
& = t+\biggl ( \sum_{s=2j}^{\infty}\bigl (
\sum_{k=j}^{s-j}(u_{q(k+e)+1}v_{q(s-k)+1}-v_{qk+1}u_{q(s+e-k)+1})\bigr
)\biggr )t^{qr+1} \\
&+ \sum_{s=qj+j}^{\infty} \biggl ( \sum \Sb m,k\geq j, \\ qm+k=s \endSb
{-(k+e)u_{q(k+e)+1}v_{qm+1}} \biggr ) t^{q(s+e)+1} + \dots
\endsplit
$$
where $r=s+e$ in the first sum. Notice now that the first sum gives an
identically zero expression and so in actual fact
$$
\split
t+v\kl u-u\kl v
&= t+ \sum_{s=qj+j}^{\infty} \biggl ( \sum \Sb m,k\geq j, \\ qm+k=s \endSb
{-(k+e)u_{q(k+e)+1}v_{qm+1}} \biggr ) t^{q(s+e)+1} +  \dots \\
&=t-iu_{qi+1}v_{qj+1}t^{q(qj+i)+1}+\dots
\endsplit
$$
and the Proposition is proved.

\enddemo

\vskip 0.5in
\heading{\bf { Section 4: $T$ is Hereditarily Just Infinite}}
\endheading
\vskip 0.25in

Now that we have established some basic commutator relationships, we are in
a position to give an alternative proof to Fesenko\rq s Lemmas about the
nature of commutators of particular elements in $T.$ For future reference,
we include a combinatorial Lemma that is proved in \cite {F1} and will be
useful to us. We use the following notation: $j=j^{\p}p^{n(j)}$ where
$j^{\p}$ is coprime to $j.$ Recall also that $q=p^r.$

\proclaim {Lemma 4.1 \cite{F1, Lemma 1}}
Fix $s$ so that $1\leq s\leq r.$ Let $i>j \geq q^2$ and let $i$ be coprime
to $p.$ Let $i_m,j_m$ satisfy the following conditions:
\roster
\item $i_m\geq i, j_m \geq j-q;$
\item $(i_m,i) =1;$
\item $j_m\geq j$ if $i_m=i;$ $qj_m+p^si_m>qj+p^si$ if $i_m>i;$
\item if $ i_m+qj_m < j+qi,$ then $i_m=r_mi+s_mq$ for integers $r_m \geq 1,
s_m \geq 0;$
\endroster
Let $v_m,w_m,x_m,y_m,z_m$ be non-negative integers so that $v_m >0 \iff
w_m>0, x_m>0 \iff y_m>0$ and $z_m >0$ only if $x_m>0.$ Let $q$ divide $z_m$
if $z_m \geq qj.$

Then the equality
$$
\split
\sum (v_mi_m+w_mqj_m)
&+\sum (x_mj_m+y_mqi_mp^{n(j_m)} +z_m) \\
& = I+qj, \quad p^{s-1}i<I\leq p^sI, p^s|I
\endsplit
$$
implies that
$$
I=p^si;
$$

Furthermore:
\newline
if $p^s<q$ then up to renumbering we have $v_1=p^s,w_1=1,i_1=i,j_1=j$ and
$v_m=w_m$ for $m>1, x_m=y_m=z_m =0$ for $m \geq 1;$
\newline
if $p^s=q$ then either up to renumbering $v_1=q,w_1=1i_1=i,j_1=j$ and
everything else is zero, or up to renumbering $x_1=q,y_1=1,i_1=i,j_1=j$ and
$v_m=w_m=z_m=0$ for $m\geq 1, x_m=y_m =0$ for $m>1.$

\endproclaim

The relevance of this result will become apparent in due course.

Fesenko considers elements $u,v \in T$ whose coefficients satisfy particular
arithmetic requirements. He takes an element $v \in T_j \backslash T_{j+1}$
for some $j \geq q^2$ as
$$
v=t+\sum_{k\geq j}v_{qk+1}t^{qk+1}
$$
where $v_{qk+1}=0$ if $j+1 \leq k \leq qj$ is not divisible by $q.$ Also
$$
u=t+u_{qi+1}t^{qi+1}
$$
for some non-zero $u_{qi+1}$ and $i>j$ relatively prime to $p.$ He proves

\proclaim{Lemma 4.2 \cite{F1, Lemma 3}} With the same notation as above
$[v,u]$ is congruent modulo $t^{1+q^2(i+j)+q}$ to
\roster
\item $t+\sum_{j_m\geq j}c_mt^{1+q(v_mi+w_mqj_m)} + \sum_{j_m\geq
j}d_mt^{1+q(x_mj_m+y_mqip^{n(j_m)}+z_m)};$
\item $t+\sum_{v\geq i} e_vt^{1+q(qj+v)}. $
\endroster
In (1) the $v_m, \dots ,z_m$ together with $i_m=i, j_m\geq j$ satisfy the
conditions of
\newline Lemma 4.1. In (2) the $e_v$ satisfy
\roster
\item if $v+qj<j+qi$ and $e_v \not=0,$ then $v=s_vi+r_vq$ for  $s_v\geq 1,
r_v \geq 0;$
\item $e_{p^si} =-iu_{qi+1}v_{qj+1}$ for $ 0 \leq s<r,$ and
$e_{qi}=(j-i)u_{qi+1}v_{qj+1}.$
\endroster
\endproclaim
Using the methods we have developed above we are now able to give an
alternative proof of this result. We prove the results about the nature of
the $e_v,$ the other results follow in a similar way.

\demo{Proof}
We demonstrate that the second description of the commutator is correct.
Essentially to prove this result we merely solve recurrence relations on the
coefficients of the commutator, which we evaluate using the methods outlined
previously. In this special case that $u,v $ have a particularly simple form
$v\kl u, u\kl v$ become easier to describe when we work modulo
$t^{q(qj+qi)+q+1}$:
$$
\split
v\kl u &=v \kl (t+u_it^{qi+1}) \\
&=t+u_it^{qi+1} + \sum_{k\geq j}v_{qk+1}(t+u_it^{qi+1})^{qk+1} \\
&= t+u_it^{qi+1}+\sum_{k\geq j}v_{qk+1}t^{qk+1} + \sum_{k\geq j}
v_{qk+1}u_it^{q(k+i)+1} \\
& +\binom {qj+1}q v_{qj+1}u_i^qt^{q(qi+j)+1} + \binom
{qj+1}{q+1}v_{qj+1}u_i^{q+1}t^{q(qi+i+j)+1} \\
&+ \binom {qj+1}{2q} v_{qj+1}u_i^{2q}t^{q(2qi+j)+1} + \dots
\endsplit
$$
and
$$
\split
u\kl v &= u\kl \bigl (t+\sum_{k \geq j} v_{qk+1}t^{qk+1} \bigr ) \\
&= t+\sum_{k\geq j} v_{qk+1}t^{qk+1} +u_i\bigl (t+\sum_{k\geq j}
v_{qk+1}t^{qk+1} \bigr )^{qi+1} \\
&= t+\sum_{k\geq j} v_{qk+1}t^{qk+1} +u_i\biggl (t^{qi+1} +\sum_{k\geq j}
v_{qk+1}t^{q(k+i)+1} \\
&+ \sum_{k\geq j} \binom {qi+1}q v_{qk+1}^qt^{q(qk+i)+1}\\ & + \sum \Sb
k\geq j \\ l \geq j \endSb \binom {qi+1}q \binom {q(i-1)+1}1
v_{qk+1}v_{ql+1}t^{q(qk+i+l)+1}
+\dots \biggr ).
\endsplit
$$
So as before we have $v\kl u-u\kl v = \sum \alpha_n(u\kl v)^{qn+1}$ where
$[v,u]=t+\sum \alpha_nt^{qn+1}$ and so
$$
\split
& \binom {qj+1}q v_{qj+1}u_i^qt^{q(qi+j)+1} + \binom
{qj+1}{q+1}v_{qj+1}u_i^{q+1}t^{q(qi+i+j)+1} \\
&+ \binom {qj+1}{2q} v_{qj+1}u_i^{2q}t^{q(2qi+j)+1} + \dots
-u_i \biggl (\sum_{k\geq j} \binom {qi+1}q v_{qk+1}^qt^{q(qk+i)+1}\\
& + \sum \Sb k\geq j \\ l \geq j \endSb \binom {qi+1}q \binom {q(i-1)+1}1
v_{qk+1}v_{ql+1}t^{q(qk+i+l)+1} +\dots \biggr ) \\
&= \sum_{n\geq K}\alpha_n \biggl ( t+\sum_{k\geq j}v_{qk+1}t^{qk+1} +u_i
\bigl ( t+\sum_{k\geq j} v_{qk+1}t^{qk+1} \bigr )^{qi+1} \biggr )^{qn+1}.
\endsplit
$$
In order to prove the first claim of the Lemma we work modulo
$t^{q(qi+qj)+1}$ and so this expression simplifies to
$$
\split
& -u_i \biggl (\sum_{k\geq j} \binom {qi+1}q v_{qk+1}^qt^{q(qk+i)+1}\\
& + \sum \Sb k\geq j \\ l \geq j \endSb \binom {qi+1}q \binom {q(i-1)+1}1
v_{qk+1}v_{ql+1}t^{q(qk+i+l)+1} +\dots \biggr ) \\
&= \sum_{n\geq K}\alpha_n \biggl ( t+\sum_{k\geq j}v_{qk+1}t^{qk+1} +u_i
\bigl ( t+\sum_{k\geq j} v_{qk+1}t^{qk+1} \bigr )^{qi+1} \biggr )^{qn+1}.
\endsplit
$$
Now we merely equate coefficients. We can trivially confirm what we already
know, namely that $\alpha_n =0$ for all $n <qj+i$ and that
$\alpha_{qj+i}=-iu_iv_{qj+1}.$
It is also immediate that $\alpha_n=0$ for all $qj+i <n<qj+i+q.$ Comparing
coefficients of  $t^{q(qj+i+q)+1}$ we see that
$$
-iu_iv_{q(j+1)+1}=\alpha_{qj+i+q};
$$
this process continues until we reach $t^{q(qj+i+j)+1}$. Here we have
$$
-iu_iv_{qj+1}^2=\alpha_{qj+i+j}+v_{qj+1}\alpha_{qj+i}
$$
and so $\alpha_{qj+i+j}=0$ as required.

At this stage it is worth tidying up what we know. The previous large
expression now simplifies to:
$$
\split
& -u_i \biggl (\sum_{k\geq j} \binom {qi+1}q v_{qk+1}^qt^{q(qk+i)+1}\\
& + \sum \Sb k\geq j \\ l \geq j \endSb \binom {qi+1}q \binom {q(i-1)+1}1
v_{qk+1}v_{ql+1}t^{q(qk+i+l)+1} +\dots \biggr ) \\
& = \alpha_{qj+i} \biggl ( t+\sum_{k\geq j}v_{qk+1}t^{qk+1} +u_i \bigl (
t+\sum_{k\geq j} v_{qk+1}t^{qk+1} \bigr )^{qi+1} \biggr )^{q(qj+i)+1} \\
&+ \alpha_{qj+i+q} \biggl ( t+\sum_{k\geq j}v_{qk+1}t^{qk+1} +\dots \biggr
)^{q(qj+i+q)+1} \\
&+ \alpha_{qj+i+2q} \biggl ( t+\sum_{k\geq j}v_{qk+1}t^{qk+1} +\dots \biggr
)^{q(qj+i+2q)+1} \\
&+\dots
\endsplit
$$
Compare coefficients of $t^{q(qj+2i)+1}:$
$$
0=\alpha_{qj+2i}+u_i\alpha_{qj+i}
$$
and so $\alpha_{qj+2i}=iu_i^2v_{qj+1}.$

Now using the nature of the non-zero coefficients of $v$ and continuing in
this way we may deduce that $\alpha_{\Delta} \not= 0 \implies \Delta =
qj+s_vi+r_vq$ for some $s_v \geq 1, r_v \geq 0$ whenever $\Delta <qi+j.$ We
can also see that for all $p^s<q,$ we have
$$
\alpha_{qj+p^si} =-iu_i^{p^s}v_{qj+1} = -iu_iv_{qj+1}.
$$
We finally need to look at $\alpha_{q(i+j)}.$ Notice that
$$
\alpha_{qi+j} =ju_iv_{qj+1}
$$
and so the same reasoning as above tells us that
$$
\alpha_{qi+p^sj}=ju_iv_{qj+1}^{p^s}=ju_iv_{qj+1}.
$$
The value of $\alpha_{q(i+j)}$ follows and the result is proved.

\enddemo

The following Lemma, amended from \cite {F1} is the final result needed in
order to demonstrate  the main thrust of Fesenko\rq s reasoning.

\proclaim{Lemma 4.3 \cite {F1, Lemma 5}}
Let $1\not= H \triangleleft_c U<_o T.$ Then for all sufficiently large $j$
coprime to $p,$ for all non-zero $a \in \Bbb F_p,$ there exists a series
$$
t+\sum_{k\geq j} a_kt^{qk+1} \in H, \quad a_j=a
$$
so that for $j+1\leq k \leq qj+q^2,$ $a_k=0$ if $q$ does not divide $k.$
\endproclaim

In the interests of completness we
now outline how Fesenko completed his proof that $T$ is hereditarily just
infinite.

\proclaim{Theorem 4.4 \cite{F1}}
$T$ is a hereditarily just infinite pro-$p$ group for $p>2$.
\endproclaim

\demo{Proof Outline}
Let $i,j,i-j$ be coprime to $p.$ Then by Lemma 4.3, given any closed
normal subgroup $H$  of an open subgroup $U$ of $T$ there exists an element
$v$ of the form described in Lemma 4.2. So applying Lemma 4.2 twice,
once to $(j,i)$ and once to $(j-p,i+qp)$ gives us modulo a high power of $t$
elements
$$
t+\sum_{v\geq i}e_vt^{1+q(qj+v)}, \quad t+\sum_{v\geq i+qp}
f_vt^{1+q(qj+v-qp)}
$$
in $H$ for arbitrary non-zero elements $e_{i},f_{i+qp}.$ Thus we may pick
these elements in such a way that the composition gives us, again modulo a
high power of $t$
$$
[v(j),u(i)] \kl [v(j-p),u(i+qp)] = t+ \sum_{v>i} g_vt^{1+q(qj+v)}.
$$
Fesenko shows that the coefficients of this power series also satisfies the
conditions of Lemma 4.2. Using this fact one can continue in this way to
produce an element
$$
t+\sum_{v\geq pi}h_vt^{1+q(qj+v)}
$$ in $H$ where again $h_{pi} \not=0$ and the $h_v$ satisfy the conditions
of Lemma 4.2. By induction, one may now produce an element of the form
$t+t^{1+q(qj+qi)}+ \dots .$ Combined with Proposition 3.4 it follows that
we may, for all $\la$ sufficiently large produce the element
$t+t^{q\la+1}+\dots $ in $H.$

Thus it follows that $T$ is HJI as claimed.
\enddemo

{\bf Remark.} The stumbling block to a proof that $T$ is hereditarily just 
infinite for
$p=2$ is the first arithmetic condition that $i,j,i-j$ are all coprime to
$p.$ Lemmas 4.1, 4.2, 4.3 remain valid as they stand for $p=2.$ It is 
entirely likely that one will be able to remove this arithmetic condition on 
$i,j$ and prove $T$ is HJI for $p=2$ without adopting that different a 
method to the one outlined above. At the moment this remains out of reach 
but I feel sure a closer investigation of the methods used in this chapter 
would bear fruit.

\vskip 0.5in
\heading{\bf {Section 5: Commutators in $S$}}
\endheading
\vskip 0.25in

$S$ is the group consisting of formal power series of the form
$$
u(t) = t+a_1t^p+a_2t^{p+1}+ \dots .
$$
We may define a filtration on $S=S_1>S_2> \dots$ where now
$$
S_{2n} := \left\{t+a_{np+1}t^{np+1}+\dots :a_i \in \Bbb F_p \right\}
$$
and
$$
S_{2n-1}:=\{t+a_{np}t^{np}+\dots :a_i \in \Bbb F_p\}.
$$

Suppose without loss of generality we are given two elements
$u=t+u_{pi}t^{pi}+\dots $ and $v=t+v_{pj}t^{pj} +\dots $ with $i=j+e.$ Then
as above one can easily verify that
$$
\split
u\kl v (t)
&= u+v-t + \sum_{s\geq i-1+jp}\biggl( \sum_{k\geq j+e}u_{kp}\sum_{j(1)
+\dots +j(kp)=sp}v_{j(1)} \dots v_{j(kp)} \biggr) t^{sp} \\
&+ \sum_{\Delta \geq (i+j)p}\biggl(\sum_{k \geq j+e}u_{kp+1} \sum _{j_(1)+
\dots +j(kp+1) =\Delta} v_{j(1)}\dots v_{j(kp+1)} \biggr)t^{\Delta}
\endsplit
$$
and similarly
$$
\split
v\kl u (t)
&= u+v-t + \sum_{s\geq j-1+ip}\biggl( \sum_{k\geq j}v_{kp}\sum_{j(1) +\dots
+j(kp)=sp}u_{j(1)} \dots u_{j(kp)} \biggr) t^{sp} \\
&+ \sum_{\Delta \geq (i+j)p}\biggl(\sum_{k \geq j}v_{kp+1} \sum _{j_(1)+
\dots +j(kp+1) =\Delta} u_{j(1)}\dots u_{j(kp+1)} \biggr)t^{\Delta}.
\endsplit
$$
\newline
In the first of these two expressions $j(l) >1 \implies j(l) \geq jp$ and
$j(l) \equiv   0,1   \mod p$ whereas in the second of the two expressions
$j(l)>1 \implies j(l) \geq ip.$

In order to prove that $S$ is HJI we must evaluate some commutators. The
following  Lemmas are simple to verify using the same method as before.

\proclaim{Lemma 5.1}Let $i>j.$
\roster
\item Let $u=t+u_it^{pi+1}, v=t+v_jt^{pj+1} +\dots .$ Then $[v,u]
=t-iu_iv_jt^{p(pj+i)+1}+\dots ;$
\item Let $u=t+u_it^{pi+1}, v=t+v_jt^{pj}+\dots .$ Then $[v,u]
=t-u_iv_jt^{p(i+j)}+ \dots ;$
\item Let $u=t+u_it^{pi}, v=t+v_jt^{pj+1}+\dots .$ Then 
$[v,u]=t+u_iv_jt^{p(i+j)}+
\dots .$
\endroster
\endproclaim

This Lemma is proved exactly as Proposition 3.4 was proved and so we omit
the details. The following is less obvious and more important to the proof
of the main result.

\proclaim{Lemma 5.2} Let $i>j$ and let $i$ be coprime to $p.$ Take
elements $u=t+u_it^{pi+1}, v=t+v_{pj}t^{pj}+ \dots $ with $u_iv_{pj}\not=0.$
Then the first power of $t$ in $[v,u]$ that has non-zero coefficient and is
congruent to $1$ modulo $p$ is
$$
-iu_iv_{pj}t^{p(pj+i-1)+1}.
$$
\endproclaim

\demo{Proof}
We proceed exactly as we have done previously. A simple exercise in
combinatorics enables us to evaluate
$$
\split
u\kl v
&= t+v_{pj}t^{pj}+v_{pj+1}t^{pj+1}+\dots
+ u(  t+v_{pj}t^{pj}+v_{pj+1}t^{pj+1}+\dots ) \\
& = u+v-t +\sum \Sb \la \geq pj \\ \la \equiv   0,1   \mod p \endSb
u_iv_{\la}t^{pi+\la}
+ \sum \Sb \la \geq pj \\ \la \equiv   0,1   modp \endSb {\binom {pi+1}p
u_iv_{\la}t^{p\la +pi+1-p}  } \\
&+ \dots
\endsplit
$$
and similarly
$$
\split
v\kl u
&= (t+u_it^{pi+1})+v_{pj}(t+u_it^{pi+1})^{pj} +\dots \\ & = u+v-t +\sum \Sb
\mu \geq pj \\ \mu \equiv   1  \mod p \endSb u_iv_{\la}t^{pi+\mu}
+ \sum \Sb \mu \geq pj \\ \mu \equiv   0,1   \mod p \endSb {\binom {\mu}p
u_iv_{\mu}t^{p^2i+\mu-p}} \\
&+ \dots
\endsplit
$$
Thus it follows that
$$
\split
v\kl u -u\kl v
&= - \sum \Sb \la \geq pj \\ \la \equiv   0  \mod p \endSb
u_iv_{\la}t^{pi+\la}
-iu_iv_{pj}t^{p(pj+i-1)+1} +\dots \\
&:= \sum_{k \in \bbbn} \alpha_kt^k
\endsplit
$$
where again $[v,u]=t+\sum_{k >1}\alpha_kt^k.$
\newline
The result is now self evident, once it is appreciated that $\alpha_k =0$
unless we have that $k$ is $0$ or $1$ modulo $p.$
\enddemo

In order to prove the main results, the last piece of the jigsaw we need is
a description of some of the torsion in $B.$ The following Lemma is simple
and indeed is contained in \cite {C2}:

\proclaim{Lemma 5.3}
Write $v^p$ for group composition of $v$ with itself $p$ times.
\roster
\item
Let $v=t+v_1t^{pn}+\dots \in S_{2n-1}\backslash S_{2n}.$ Then
$$
v^p \in S_{2np-1};
$$
\item
Let $v=t+v_1t^{pn+1}+\dots \in S_{2n}\backslash S_{2n+1}.$ Then
$$
v^p \in S_{2np}.
$$
\endroster
\endproclaim

In his thesis \cite{Y} York gives a complete description of the elements of
order $p$ in the Nottingham Group. The next two results are those contained
therein that are of relevance to the situation here.

\proclaim {Theorem 5.4, \cite {Y, Theorem 5.5.3}} Let $\alpha =t+\sum_{k
\geq n}a_{pk}t^{pk} \in \Cal J(\Bbb F_p).$ Then $\alpha$ has order $p$ if
and only if
$$
a_{(2n+s-np+1)p-1}=f_s(a_{(2n+s-np+1)p-2},\dots , a_{np})
$$
for some given polynomial $f_s$ dependent upon $s, \alpha $ and $n+s \geq
np.$
\endproclaim

\proclaim{Theorem 5.5 \cite{Y, Theorem 5.5.4}}
Let $\alpha = t+at^{pk+1}+\dots \in \Cal J(\Bbb F_p).$ Then $\alpha$ has
infinite order and $\alpha^p = t+at^{p(pk)+1}+\dots .$
\endproclaim
These results follow merely from a close examination of the coefficients of
the power series.
Notice also that Theorem 5.5 strengthens Lemma 5.3(2) above to show that
given $v \in S_{2n}\backslash S_{2n+1}$ then $v^p \in S_{2pn}\backslash
S_{2pn+1}.$

Given these results we are in a position to prove the final Proposition we
will need in order to prove that $B$ is hereditarily just infinite.

\proclaim{Proposition 5.6}
Let $1 \not= H \triangleleft _c U<_oB.$ Then $H$ is infinite and
furthermore, $H$ contains elements of arbitrarily large depth.
\endproclaim

\demo{Proof}
As $H$ is non-trivial there is an element in $H$ different from $t.$ Take
such an element $v.$ If $v=t+v_jt^{pj+1}+\dots $ then by Theorem 5.5 $v$
has
infinite order and $H$ is infinite. Thus taking powers of $v$ gives elements
of arbitrary depth in $H.$ Thus in this case we\rq re done.

Suppose instead that $v=t+v_jt^{pj}+\dots .$ Then Lemma 5.1(2) implies
that for any $i>j,$
$$
[t+ut^{pi+1},v]=t+uv_jt^{p(i+j)}+\dots
$$
and so it follows that $H$ is infinite.

The same commutator relation tells us that we may alter coefficients of
powers of $t$ in $v$ occuring after the $p(2j+1)$-th power of $t$ by
composing $v$ with
$$
[t+u_{j}t^{pj+1},v] =t+u_jv_jt^{p(2j+1)}+\dots.
$$
Thus by Theorem 5.4, $H$ contains an element not of order $p$ where the
coefficients of later powers of $t$ are determined by the previous ones.
Continue this process with the $p$-th power of this element. The result
follows.

\enddemo

\vskip 0.5in
\heading{{\bf Section 6: $S$ is Hereditarily Just Infinite}}
\endheading
\vskip 0.25in

The proof mirrors Fesenko\rq s proof outlined earlier. In particular, we
take a non-trivial (closed)  normal subgroup $H$  of an open subgroup $U$ of
$S.$ Then $H$ is infinite and contains elements of arbitrary depth. We take
such an element $v$ and show that by taking appropriate commutators of $v,$
and of powers of $v$ with arbitrary elements  $u \in S_i$ for sufficiently
large $i$  that we may realise, for $\la$ sufficiently large in $p\bbbn \cup
p\bbbn -1$, any element of the form $t+t^{\la}+\dots $ in $H.$ This will be
sufficient to complete the proof.

\proclaim{Proposition 6.1} Suppose that $H$ is a non-trivial closed normal
subgroup of an open subgroup $U$ of $S.$ Suppose also that $H$ contains an
element of infinite order of the form
$$
v=t+v_1t^{pj+1}+ \dots
$$
for some non-zero $v_1 \in \Bbb F_p.$ Then $H$ is open.
\endproclaim

\demo{Proof}
If $v=t+v_1t^{pj+1}+\dots $ has infinite order, then taking a sufficiently
large power of $v$ we may assume that $j$ is arbitrarily large. Then by
Lemma 5.1(3), we can commutate $v$ in such a way that $v$ approximates an
element of the group $T$ arbitrarily closely. $T$ is hereditarily just
infinite, and so it follows that for all sufficiently large $\la$ we may
realise $t+t^{p\la +1}+\dots $ as an element of $H.$

Also by Lemma 5.1(3) we may realise $t+t^{p\la}+\dots $ as an element of
$H$ for all sufficiently large $\la.$ Thus the result follows.

\enddemo

In order to prove that $S$ is hereditarily just infinite, it is now
sufficient to establish the next

\proclaim{Proposition 6.2}
Let $H$ be a non-trivial normal closed subgroup of an open subgroup $U$ of
$S.$ Then there exists in $H$ an element of the form
$$
v=t+v_1t^{pj+1}+\dots , \quad v_1 \not= 0
$$
necessarily of infinite order.
\endproclaim

\demo{Proof}
Suppose that this is not the case. By the previous results we have proved,
$H$ must contain an element of the form
$$
v=t+v_jt^{pj}+\dots
$$
not of order $p.$
Then from Lemma 5.2 we have
$$
[v,t+u^{\p}t^{pi+1}]=t+u^{\p}vt^{p(i+j)} + \dots +(-iu^{\p}v)t^{p(pj+i-1)+1}
+\dots
$$
where $t^{p(pj+i-1)+1}$ is the first power of $t$ that is $1  \mod p$ and
has non-zero coefficient.

Also by Lemma 5.3 we may deduce that
$$
v^p=t +\alpha t^{p\la}+\dots
$$
for some non-zero $\alpha$ and sufficiently large $\la.$  Thus
$$
[v^p,t+ut^{pi+1}] = t+\alpha ut^{p(\la+i)} + \dots
+(-iu\alpha)t^{p(p\la+i-1)+1} +\dots
$$
Similarly
$$
[v,t+u^{\p}t^{p(\la+i-j-1)+1}] =t+vu^{\p}t^{p(\la+i)} + \dots
+(j-i)vu^{\p}t^{p(pj+\la+i-j-1)+1} +\dots
$$
where now $t^{p(pj+\la+i-j-1)+1}$ is the first power of $t$ that is $1
\mod p$ and has non-zero coefficient.  Thus noticing   taking the 
composition
of these two commutators for an appropriate choice of $u,u^{\p}$ gives us
the element
$$
t+\gamma t^{p(\la+i+1)} + \dots + \delta t^{p(pj+\la+i-j-1)+1} +\dots
$$
where $\delta$ is non-zero, and this is the coefficient of the lowest power
of $t$ that is $1  \mod p.$

We can repeat this process as often as we need, arriving in a finite number
of steps at the element
$$
t+\delta t^{p(pj+\la+i-j-1)+1}+\dots
$$
The result follows.
\enddemo

\proclaim{Corollary 6.3} $S$ is a hereditarily just infinite subgroup of
$\Cal J.$
\endproclaim

\demo{Proof}
Simply combine Propositions 6.1, 6.2.
\enddemo

{\bf Remark.} Notice that this result seems to display the characteristics 
that one would
expect in attempting to prove a Theorem of this nature. It is generally
harder to produce as a commutator elements of the form $t+t^{pi+1}+\dots$
than anything else. It was this problem that caused difficulties when people
tried to prove the Nottingham Group is HJI for $p=2.$ Notice that when $p=2$
we actually have that $\Cal J =S.$ These difficulties have been overcome for
$\Cal J$ in a recent article by Hegedus \cite{H}, although it is worth
pointing out that the first proof that $\Cal J$ is HJI for $p=2$ was
communicated by Fesenko to Leedham-Green in 1999 and in fact the
result follows from reasoning in \cite {F1}. In some sense Fesenko\rq s
calculations do the hard work for us in this case, and it is merely
necessary to piece all the various fragments together as we have done here.

It was only to ease difficulties of notation, and to keep the calculations
as simple as possible that we took $B= p\bbbn \cup p\bbbn -1.$ The same
calculations will lead with no great complications to a proof that
$S_{m,n}:=\Cal J(B_{m,n})$ is hereditarily just infinite, for
$B_{m,n}:=p^n\bbbn \cup p^m\bbbn -1.$ However in the light of the method of
proof, this reasoning is only valid for odd primes as we do not as yet have
a full even characteristic equivalent to Theorem 4.4.

Lastly note that the group $S$ has the property that every open subgroup has 
non-trivial torsion. Thus the group $S$ remains as a potential candidate to 
be the Galois group of a just infinite unramified field extension of $\Bbb 
Q.$ See \cite {B} for more details.

\vskip 0.5in
\heading{\bf {Appendix}}
\endheading
\vskip 0.25in

Since the completion of this work, a preprint of Barnea and Klopsch
\cite{BK} has appeared which also considers various properties of index
subgroups in $\Cal J. $ The main results contained in this paper are
summarised in this

\proclaim {Theorem 7 {\cite{BK}}}
Let $S_{r,s} := \Cal J (p^r\bbbn \cup p^s\bbbn -1).$ Then
\roster
\item $S_{r,s}$ is hereditarily just infinite$;$
\item $S_{r,r}$ is of finite width and infinite obliquity$;$
\item The width of $S_{r,r}$ is $\leq p+p^r-1;$
\item The groups $S_{r,r}$ are pairwise non-commensurable.
\endroster
\endproclaim

They also consider implications for the Hausdorff Spectrum of the Nottingham
Group, and indeed calculate a large part of this spectrum.

To prove results about the normal subgroup structure of a substitution group
of formal power series, the only real possible approach to this problem is
to evaluate commutators. In \cite{BK} the authors adopt a different method
to do this than the one contained here, and it is reassuring to see that our
conclusions are the same.

In the light of \cite{BK} a major problem to tackle would still seem to be
to calculate the width of the Fesenko Group $T.$ It does not appear possible
to use the calculations of Barnea and Klopsch for the width of $S_{r,r}$ to
do this. Indeed it is still far from clear to me whether or not $T$ will
have finite width. Fesenko believes that $T$ will have finite width based on
some lengthy preliminary calculations he performed several years ago but he
emphasizes this has still to be
confirmed. However that $S$ has infinite obliquity would certainly suggest
that $T$ will do too.

\vskip 0.5in
\heading{Bibliography}
\endheading

\
\newline
[B] N Boston, {\it Some cases of the Fontaine-Mazur Conjecture. II,} Alg.
Number
Th. Arch. {\bf 129} (1998).
\newline
[BK] Y Barnea, B Klopsch, {\it Index subgroups of the Nottingham Group,}
Preprint, (2002).
\newline
[C1] R Camina, {\it Subgroups of the Nottingham Group,} J. Algebra {\bf 196}
(1997), 101--113.
\newline
[C2] R Camina, {\it The Nottingham Group,} New Horizons in Pro-$p$ Groups,
\linebreak Birkhauser Publ., 2000.
\newline
[CG] J Coates, R Greenberg, {\it Kummer theory for the abelian varieties 
over local fields,} Invent. Math. {\bf 124} (1996), 129--174.
\newline
[duSF] M du Sautoy, I Fesenko, {\it Where the wild things are:
ramification
groups and the Nottingham Group,} New Horizons in Pro-$p$ Groups, Birkhauser
Publ., 2000, 285--326.
\newline
[F1] I Fesenko, {\it On just infinite pro-p groups and arithmetically
profinite
extensions of local fields,} J. reine angew. Math. {\bf 517} (1999), 61--80.
\newline
[F2] I Fesenko, {\it On deeply ramified extensions,} J. LMS (2) {\bf 57}
(1998),
325--335.
\newline
[FV] I Fesenko, S Vostokov, {\it Local fields and their extensions,} AMS,
Providence, R.I. 1993.
\newline
[H] P Hegedus, {\it The Nottingham group for $p=2$,} preprint, 2001.
\newline
[J] S A Jennings, {\it Substitution groups of formal power series under
substitution,} Canad. J. Math. {\bf 6} (1954), 325--340.
\newline
[Jo] D Johnson, {\it The group of formal power series under substitution,}
Austral. Math. Soc. {\bf 45} (1988), 296--302.
\newline
[KLP] G Klaas, C R Leedham-Green, W Plesken, {\it Linear pro-p groups of
finite
width,} preprint, Aachen-London, 1997.
\newline
[K] B Klopsch, {\it Substitution groups, subgroup growth and other topics,}
D.Phil Thesis, Oxford, 1999.
\newline
[Lu] E Lucas, {\it Sur les congruences des nombres eulerians et des
coefficients differentiels des fonctions trigonometriques, suivant un module
premier,} Bull. Soc. Math. France {\bf 6} (1878), 49--54.
\newline
[Y] I York {\it The group of formal power series under substitution,} PhD
Thesis, Nottingham, 1990.

\enddocument